\documentclass[a4paper,twopage,reqno,12pt]{amsart}
\usepackage[top=30mm,right=30mm,bottom=30mm,left=30mm]{geometry}

\usepackage{tabularx}
\usepackage{graphicx}
\usepackage{amsmath}
\usepackage{amssymb}
\usepackage{amsfonts}
\usepackage{amsthm}
\usepackage{amstext}
\usepackage{amsbsy}
\usepackage{amsopn}
\usepackage{amscd}
\usepackage{enumerate}
\usepackage{color}
\usepackage[colorlinks]{hyperref}
\usepackage[hyperpageref]{backref}
\usepackage{multirow}
\usepackage{longtable}
\usepackage{float}
\usepackage{lscape}
\usepackage{pdflscape}
\usepackage{url}
\usepackage[numbers,sort&compress]{natbib}

\newtheorem{theorem}{Theorem}[section]

\theoremstyle{definition}

\numberwithin{equation}{section}

\newcommand{\GF}{\mathrm{GF}}

\newcommand{\GU}{\mathrm{GU}}
\newcommand{\G}{\mathrm{G}}

\newcommand{\GL}{\mathrm{GL}}
\newcommand{\SL}{\mathrm{SL}}

\newcommand{\GSp}{\mathrm{GSp}}

\newcommand{\SU}{\mathrm{SU}}

\newcommand{\PSL}{\mathrm{PSL}}
\newcommand{\PSU}{\mathrm{PSU}}

\newcommand{\AGL}{\mathrm{AGL}}

\newcommand{\PGaL}{\mathrm{P\Gamma L}}
\newcommand{\GaL}{\mathrm{\Gamma L}}

\newcommand{\PGaU}{\mathrm{P\Gamma U}}

\newcommand{\AGaL}{\mathrm{A\Gamma L}}

\newcommand{\A}{\mathrm{A}}
\newcommand{\W}{\mathrm{W}}

\renewcommand{\S}{\mathrm{S}}
\newcommand{\Q}{\mathrm{Q}}

\newcommand{\D}{\mathrm{D}}

\newcommand{\Aut}{\mathrm{Aut}}

\newcommand{\PG}{\mathrm{PG}}
\newcommand{\AG}{\mathrm{AG}}

\newcommand{\B}{\mathrm{B}}

\newcommand{\M}{\mathrm{M}}

\newcommand{\Fbb}{\mathbb{F}}

\newcommand{\Dmc}{\mathcal{D}}
\newcommand{\Bmc}{\mathcal{B}}

\newcommand{\Pmc}{\mathcal{P}}

\newcommand{\Umc}{\mathcal{U}}

\renewcommand{\leq}{\leqslant}
\renewcommand{\geq}{\geqslant}

\begin{document}
\title[Block designs admitting flag-transitive automorphism groups]{Block designs with $\gcd(r,\lambda)=1$ admitting flag-transitive automorphism groups}

\author[]{
	Seyed Hassan Alavi,  
	Mohsen Bayat, 
	Mauro Biliotti, 
	Ashraf Daneshkhah, 
	Eliana Francot,  
	Haiyan Guan, 
	Alessandro Montinaro,  
	Fatemeh Mouseli, 
	Pierluigi Rizzo, 
	Delu Tian, 
	Yajie Wang, 
	Xiaoqin Zhan, 
	Yongli Zhang,   
	Shenglin Zhou and
	Yan Zhu
}

\thanks{Corresponding author: Ashraf Daneshkhah}
\address{Seyed Hassan Alavi, Mohsen Bayat, Ashraf Daneshkhah$^{\ast}$ and Fatemeh Mouseli, Department of Mathematics, Faculty of Science, Bu-Ali Sina University, Hamedan, Iran.}%
\email{alavi.s.hassan@basu.ac.ir and alavi.s.hassan@gmail.com}
\email{m.bayat@sci.basu.ac.ir}
\email{adanesh@basu.ac.ir}
\email{f.mooseli@gmail.com}

\address{Mauro Biliotti, Eliana Francot, Alessandro Montinaro and Pierluigi Rizzo, Dipartimento di Matematica e Fisica “E. De Giorgi”, University of Salento, Lecce, Italy}
\email{alessandro.montinaro@unisalento.it}
\email{mauro.biliotti@unisalento.it}
\email{eliana.francot@unisalento.it}
\email{pierluigi.rizzo@unisalento.it}

\address{Haiyan Guan, Department of Mathematics, China Three Gorges University, Yichang, Hubei 443002, PR China} \email{ghyan928@ctgu.edu.cn}

\address{Delu Tian, Department of Mathematics, Guangdong University of Education, Guangzhou, Guangdong 510303, PR China}
\email{tiandelu@gdei.edu.cn}

\address{Yajie Wang, School of Mathematics and Statistics, North China University of Water Resources and Electric Power, Zhengzhou, Henan 450046, PR China}  
\email{wangyajie8786@163.com}

\address{Xiaoqin Zhan and Suyun Ding, School of Science, East China JiaoTong University, Nanchang, Jiangxi 330013, PR China} 
\email{zhanxiaoqinshuai@126.com}

\address{Yongli Zhang, Department of Mathematics, Guangdong Polytechnic Normal University, Guangzhou, Guangdong 510665,PR China} 
\email{lilyzhang0854@163.com}
\address{Shenglin Zhou, School of Mathematics, South China University of Technology, Guangzhou, Guangdong 510640, PR China}
\email{slzhou@scut.edu.cn}
\address{Yan Zhu, School of Medical Information and Engineering, Xuzhou Medical University, Xuzhou, Jiangsu 221004, PR China}
\email{zhuyan@xzhmu.edu.cn}

\subjclass[]{05B05; 05B25; 20B25}%
\keywords{ $2$-design; automorphism group}
\date{\today}%

\begin{abstract}
 In this paper, we present a classification of $2$-designs with $\gcd(r,\lambda)=1$ admitting flag-transitive automorphism groups. If $G$ is a flag-transitive automorphism group of a non-trivial $2$-design $\Dmc$ with $\gcd(r,\lambda)=1$, then either $(\Dmc,G)$ is one of the known examples described in this paper, or $\Dmc$ has $q = p^{d}$ points with $p$ prime and $G$ is a subgroup of $\AGaL_{1}(q)$.      
\end{abstract}

\maketitle
\section{Introduction}\label{sec:intro}

Our purpose here is to announce a classification of the pairs $(\Dmc, G)$, where $\Dmc$ is a nontrivial $2$-design with $\gcd(r,\lambda)=1$ and $G$ is a group of automorphisms acting transitively on the flags of $\Dmc$. Here, a $2$-design $\Dmc$ with parameters $(v,k,\lambda)$ is a pair $(\Pmc,\Bmc)$ with a set $\Pmc$ of $v$ points and a set $\Bmc$ of $b$ blocks such that each block is a $k$-subset of $\Pmc$ and each two distinct points are contained in $\lambda$ blocks. A \emph{symmetric} design is a $2$-design with the same number of points and blocks, that is to say, $v=b$. The \emph{replication number} $r$ of $\Dmc$ is the number of blocks incident with a given point. A \emph{flag} of $\Dmc$ is a point-block pair $(\alpha, B)$ such that $\alpha \in B$. An automorphism group $G$ of $\Dmc$ is a group of permutation on $\Pmc$ which maps blocks to blocks and preserving the incidence. Further notation and definitions in both design theory and group theory are standard and can be found, for example, in~\cite{b:Atlas,b:Dixon,t:Kleidman,b:lander}. In this paper, we give a classification of $2$-designs with $\gcd(r,\lambda)=1$. 

\begin{theorem}\label{thm:main}
	Suppose that $G$ is a flag-transitive automorphism group of a non-trivial $2$-design $\Dmc$ with $\gcd(r,\lambda)=1$. Then either  
	\begin{enumerate}[\rm (a)]
		\item $(\Dmc,G)$ is one of the known examples described in {\rm Sections \ref{sec:AS} and \ref{sec:HA}} below, or
		\item $\Dmc$ has $q = p^{d}$ points with $p$ prime and $G$ is a subgroup of  $\AGaL_{1}(q)$.
	\end{enumerate}
\end{theorem}

If $\lambda=1$, then $\Dmc$ is a linear space and a result of Higman and McLaughlin \cite{a:HigMcL-61} shows that $G$ acts primitively on the points of $\Dmc$, and using the O'Nan-Scott theorem for finite primitive permutation groups,  Buekenhout, Delandtsheer and Doyen proved that $G$ is of almost simple or affine type, and in 1990, a classification of linear spaces admitting flag-transitive automorphism groups has been announced in \cite{a:BDDKLS90}. The proof of this result was published in several papers. The proof in the affine case is due to Liebeck \cite{a:Liebeck-98-Affine}. The almost simple case has been treated by several authors. Delandtsheer took the case where the simple socle is an alternating group \cite{a:Delandsheer-Lin-An-2001}. She also handled the case where the group $G$ is one of the simple groups $\PSL_{2}(q)$, $\PSL_{3}(q)$, $\PSU_{3}(q)$ and $^{2}\B_{2}(q)$ \cite{a:Delandtsheer-86}. Kleidman solved the case where the socle is an exceptional group of Lie type \cite{a:Kleidman-Exp} and he gave a proof for three of the ten families of exceptional groups and some hints for the remaining cases. The case of the sporadic groups was ruled out by Buekenhout, Delandtsheer, and Doyen \cite{a:BDD-spor} and Davies \cite{t:Davies87}. Finally, Saxl completed the proof in \cite{a:Saxl2002}, where he dealt with the remaining families of exceptional type together with the classical groups of Lie type.

In general, if the replication number of a $2$-design $\Dmc$ is coprime to $\lambda$, then Dembowski \cite[2.3.7]{b:Dembowski} proved  that flag-transitive automorphism groups $G$ of $\Dmc$ act primitively on the point set of $\Dmc$. In 1988, Zieschang~\cite{a:Zieschang-88} proved that such an automorphism group is of almost simple or affine type:
\begin{enumerate}[\rm (a)]
	\item Almost simple type: $G$ has a nonabelian simple normal subgroup $X$ such that $X\unlhd G\leq \Aut(X)$;
	\item Affine type: $G$ is a subgroup of the affine group $\A\GL_d(p)$ containing the translation group $T$. The socle $T$ of $G$ is an elementary abelian $p$-group of order $p^d$ with $d\geq1$. Moreover, $G=TG_{0}$, where the point-stabilizer $G_0$ of $G$ is an irreducible subgroup of $\GL_d(p)$.
\end{enumerate}

Since 2015, we have analysed these two possible cases and proved Theorem \ref{thm:main} in several papers. 
The proofs rely on the classification of finite simple groups and detailed knowledge concerning their subgroups and their permutation representations. We have not been able to resolve the case where $G\leq \AGaL_{1}(q)$, however, in Section~\ref{sec:onesemdim}, we present some examples in this case. \smallskip


\noindent\textbf{Proof of Theorem \ref{thm:main}.} 
The proof of Theorem~\ref{thm:main} appears in several papers. Suppose that $G$ is a flag-transitive automorphism group of a non-trivial $2$-design $\Dmc$ with $\gcd(r,\lambda)=1$. It follows from \cite[Theorem]{a:Zieschang-88} that such an automorphism group $G$ is of almost simple or affine type. If $G$ is of almost simple type, then $(G,\Dmc)$ is one of the examples described in Section \ref{sec:AS} below, namely, Example \ref{ex:proj-space}-\ref{ex:other}, see the main result in ~\cite{a:A-Exp-CP,a:A-Sz,a:ABD-Exp,a:ABD-Un-CP,a:ABD-Un-CP-cor,a:ADM-AS-CP,a:BDD-1988,a:Delandsheer-Lin-An-2001,a:Kantor-85-2-trans,a:Kantor-87-Odd,a:Saxl2002,a:Zhou-sym-sporadic,a:Zang-Un-CP-comment,a:Zhan-CP-nonsym-sprodic,a:Zhou-Exp-CP,a:Zhuo-CP-sym-alt,a:Zhou-CP-nonsym-alt,a:Saxl2002}. If $G$ is of affine type, then by \cite{a:Biliotti-CP-nonsol-HA,a:Biliotti-CP-sym-affine,a:Biliotti-CP-sol-HA,a:Liebeck-98-Affine}, the group $G$ is a subgroup of $\AGaL_{1}(q)$ or $(\Dmc,G)$ is as in Example \ref{ex:Desar-affine}-\ref{ex:tensprod} described in Section \ref{sec:HA} below. 

\section{Almost simple type}\label{sec:AS}

In this section, we provide some examples of $2$-designs with $\gcd(r,\lambda)=1$ admitting a flag-transitive automorphism almost simple group with socle $X$. The $2$-designs in Examples~\ref{ex:ree-designs} and~\ref{ex:suzuki} appear in \cite{a:A-Exp-CP,a:A-Sz,a:ABD-Exp,a:Zhou-Exp-CP} when the socle $X$ of $G$ is a finite simple exceptional group of Lie type. The remaining examples arose from studying linear spaces, $2$-transitive automorphism groups of $2$-designs and automorphism groups of $2$-designs with socle alternating groups, sporadic simple groups and  finite simple classical groups of Lie type \cite{a:ABD-Un-CP,a:ABD-Un-CP-cor,a:ADM-AS-CP,a:BDD-1988,a:Delandsheer-Lin-An-2001,a:Kantor-85-2-trans,a:Kantor-87-Odd,a:Saxl2002,a:Zhou-sym-sporadic,a:Zang-Un-CP-comment,a:Zhan-CP-nonsym-sprodic,a:Zhuo-CP-sym-alt,a:Zhou-CP-nonsym-alt}.  We note here that the examples of symmetric designs occur only in Examples~\ref{ex:proj-space} and~\ref{ex:other}. 

\subsection{Point-hyperplane designs}\label{ex:proj-space}
The point-hyperplane design of the projective space $\PG_{n-1}(q)$ with parameters  $((q^{n}-1)/(q-1),(q^{n-1}-1)/(q-1),(q^{n-2}-1)/(q-1))$ for $n\geq 3$ is a well-known example of flag-transitive symmetric design. A group $G$ with $\PSL_{n}(q) \leq G \leq \PGaL_{n}(q)$ acts flag-transitively on $\PG_{n-1}(q)$. If $n=3$, then we have the \emph{desarguesian plane} with parameters $(q^{2}+q+1,q+1,1)$ which is a \emph{projective plane}. The \emph{Fano plane} is the unique projective plane $\PG_2(2)$ with parameters $(7,3,1)$. We remark that there is one additional example with parameters $(15,7,3)$ and $G=\A_7$ which we view as a projective space, see \cite{a:BDD-1988,a:Kantor-87-Odd}.  

\subsection{Designs with projective points}

Suppose that $G$ is an almost simple group with socle $X=\PSL_n(q)$ with $n\geq 3$ and $H$ is a parabolic subgroup of type $P_{1}$. Let $B$ be a line (a $2$-dimensional subspace of the vector space  $V=\Fbb_{q}^{n}$) in $\PG_{n-1}(q)$ and $\alpha\in B$. Let also $\Pmc$ be the point set of $\PG_{n-1}(q)$ and $\Bmc=(B\setminus\{\alpha\})^{G}$. Then the incidence structure   
$\Dmc=(\Pmc,\Bmc)$ is a $2$-$((q^{n}-1)/(q-1),q,q-1)$ with $\gcd(n-1,q-1)=1$ admitting automorphism group $G$ with point-stabilizer $H$ such that $H\cap X\cong \,^{\hat{}}[q^{n-1}]{:}\SL_{n-1}(q){\cdot} (q-1)$, see \cite{a:ADM-AS-CP}. 

\subsection{Witt-Bose-Shrikhande spaces}\label{ex:witt}
This space is a $2$-design with parameters  $(2^{a-1}(2^{a}-1), 2^{a-1},1)$ which can be defined from the group $\PSL_2(q)$ with $q=2^{a}$ for $a\geq 3$ \cite{a:BDD-1988}. In this incidence structure which is denoted by $\W(q)$, the points are the dihedral subgroups of $\PSL_{2}(q)$ of order $2(q+1)$, the blocks are the involutions of $\PSL_2(q)$, and a point is incident with a block precisely when the dihedral subgroup contains the involution. An almost simple group $G$ with socle $X=\PSL_2(q)$ acts flag-transitively on Witt-Bose-Shrikhande space. Moreover, this space is not a symmetric design.

\subsection{Hermitian unitals}\label{ex:hermition}
The \emph{Hermitian unital} $\Umc_{H}(q)$ with parameters  $(q^{3}+1,q+1,1)$ is a well-known example of flag-transitive linear spaces \cite{a:Kantor-85-Homgeneous}. Let $V$ be a three-dimensional vector space over the field $\GF(q^2)$ with a non-degenerate Hermitian form. The Hermitian unital of order $q$ is an incidence structure whose points are $q^3+1$ totally isotropic $1$-spaces in $V$, the blocks are the sets of $q+1$ points lying in a non-degenerate $2$-space, and the incidence is given by inclusion. This structure is not symmetric and any group $G$ with $\PSU_3(q)\leq G\leq \PGaU_3(q)$ acts flag-transitively on the Hermitian unital design.

\subsection{Unitary designs}\label{ex:unitary}
A \emph{unitary design} $\Dmc=(\Pmc,\Bmc)$ with parameters $(q^{3} +1, q, q -1)$ can be constructed from the Hermitian unital as in \eqref{ex:hermition}. The point set of $\Dmc$ is the point set of the Hermitian unital $\Umc_{H}(q)$ and the block set $\Bmc$ is $(\ell\setminus\{\gamma\})^{G}$, where $\ell$ is a line of $\Umc_{H}(q)$ and $\gamma\in \ell$, see \cite{a:ABD-Un-CP-cor,a:Zang-Un-CP-comment}. This structure is not symmetric and a group $G$ with $\PSU_3(q)\leq G\leq \PGaU_3(q)$ is a flag-transitive automorphism group of $\Dmc$. 
It is worth noting by \cite{a:DLPB-2021-lam2} that there is a general construction method for $2$-designs from linear spaces, and the unitary designs described here can also be obtained in this way from the Hermitian unital design.

\subsection{Ree unitals}\label{ex:ree-unital}
The \emph{Ree Unital} $\Umc_{R}(q)$ is first discovered by L\"{u}neburg \cite{a:Luneburg-66}, and this examples arose from studying flag-transitive linear spaces \cite{a:Kantor-85-Homgeneous,a:Kleidman-Exp}. This design has parameters  $(q^{3}+1,q+1,1)$ with $q=3^a\geq 27$.   The points and blocks of $\Umc_{R}(q)$ are the Sylow $3$-subgroups and the involutions of $^{2}\!\G_{2}(q)$, respectively, and a point is incident with a block if the  block normalizes the  point. This incidence structure is a linear space and any group with $^{2}\!\G_{2}(q)\leq G\leq \Aut(^{2}\!\G_{2}(q))$ acts flag-transitively. This design is not symmetric. Note for $q=3$ that the Ree Unital $\Umc_{R}(3)$ is isomorphic to the  Witt-Bose-Shrikhande space $\W(8)$ as $^{2}\!\G_{2}(3)'$ is isomorphic to $\PSL_{2}(8)$.

\subsection{Ree designs}\label{ex:ree-designs}
Suppose that $G$ is an almost simple group with socle $X={}^2\!\G_{2}(q)$ for $q=3^a$ and $a\geq 3$ odd. Let $H$, $K_{1}$ and $K_{2}$ be subgroups of $G$ such that $H\cap X\cong q^{3}{:}(q-1)$, $K_{1}\cap X=q{:}(q-1)$ and $K_{2}\cap X \cong q^{2}{:}(q-1)$. The coset geometries $(X,H\cap X,K_{i}\cap X)$ give rise to the $2$-designs with parameters  $v=q^3+1$, $b=q^{3-i}(q^{3}+1)$, $r=q^{3}$,  $k=q^{i}$ and $\lambda=q^{i}-1$, for $i=1,2$,  see \cite{a:A-Exp-CP,a:Zhou-Exp-CP}. Since $G$ is $2$-transitive on the point set of this structure and $\gcd(r,\lambda) = 1$,  $X$ is flag-transitive \cite[2.3.8]{b:Dembowski}. Note that $H\cap K_{i} \cap X$ is a cyclic group of order $q-1$ and the subgroup $K_{i}\cap X$ has an orbit of length $q^{i}$ with $i=1,2$, see \cite{a:D-Ree,a:Zhou-Exp-CP}. If $\Pmc=\{1,\ldots,v\}$, then since $X$ is $2$-transitive, \cite[Proposition~4.6]{b:Beth-I} gives rise to a $2$-design $\Dmc_{i}=(\Pmc,B_{i}^{X})$ with parameters $(q^3+1,q^{i},q^{i}-1)$, for $i=1,2$, which is not symmetric, and the group $G$ is flag-transitive on $\Dmc_i$. An explicit construction for each of these designs is given in \cite{a:D-Ree}. For $q=27$,  in \cite[Table 1]{a:A-Exp-CP}, we introduced base blocks for these type of designs.

\subsection{Suzuki-Tits ovoid designs}\label{ex:suzuki}
These designs  arose from studying block designs with flag-transitive almost simple finite exceptional groups, see \cite{a:A-Exp-CP,a:Zhou-Exp-CP}.  
An \emph{ovoid} in $\PG_3(q)$ with $q>2$, is a set of $q^{2} + 1$ points such that no three of which are collinear. If $q$ is odd, then all ovoids are elliptic quadratics, see \cite{a:Barlotti-55,a:Panella-55}, while in even characteristic, there is only one known family of ovoids that are not elliptic quadrics in which $q\geq 8$ is an odd power of $2$. These were discovered by Tits \cite{a:Tits-ovoid}, and are now called the \emph{Suzuki-Tits ovoids} since the Suzuki groups $X={}^2\!\B_{2}(q)$ naturally act on these ovoids. A Suzuki-Tits ovoid design is a $2$-design with parameters  $(q^2+1,q,q-1)$ which is not symmetric, see \cite{a:A-Exp-CP,a:A-Sz,a:Zhou-Exp-CP}. The point set of this design is the Suzuki-Tits ovoid and the block set is the $X$-orbit $B^{X}$, where $B$ is the unique orbit of the subgroup $K:=q{:}(q-1)$ of length $q$, see \cite{a:A-Sz,a:Zhou-Exp-CP}. For $q\in \{8,32\}$, we construct these type of designs with explicit base blocks in \cite[Table~1]{a:A-Exp-CP}. There is another classical construction for these designs in geometry, by taking points as ovoids in $\PG_3(q)$ and blocks as pointed conics minus the distinguished points. 

\subsection{More examples}\label{ex:other}
Table~\ref{tbl:main} illustrates eleven examples of designs. Each design $\Dmc=(\Pmc, \Bmc)$ with parameters $(v,k,\lambda)$ is the unique design with flag-transitive automorphism
group $G$ as in the seventh column of Table~\ref{tbl:main}. The point-stabilizer and block-stabilizer of $\Dmc$ are also given in the same table with appropriate references in each case, see also \cite{a:ADM-AS-CP} for explicit base blocks. 
\begin{table}
	\caption{ Some nontrivial $2$-designs $\Dmc$ with $\gcd(r,\lambda)=1$ admitting flag-transitive and point-primitive automorphism group $G$.}\label{tbl:main}
	\resizebox{\textwidth}{!}{
		\begin{tabular}{clllllllllll}
			\noalign{\smallskip}\hline\noalign{\smallskip}
			Line &
			$v$ &
			$b$ &
			$r$ &
			$k$ &
			$\lambda$ &
			$G$ &
			$G_{\alpha}$&
			$G_{B}$ &
			$\Aut(\Dmc)$ &
			Design & 
			References
			\\
			\noalign{\smallskip}\hline\noalign{\smallskip}
			$1$ &
			$6$ &
			$10$ &
			$5$ &
			$3$ &
			$2$ &
			$\PSL_{2}(5)$ &
			$\D_{10}$ &
			$\S_{3}$ &
			$\PSL_{2}(5)$ &
			&
			\cite{b:Handbook,a:Zhou-PSL2-non-sym,a:Zhou-CP-nonsym-alt}
			\\
			$2$ &
			$8$ &
			$14$ &
			$7$ &
			$4$ &
			$3$ &
			$\PSL_{2}(7)$&
			$7{:}3$ &
			$\A_{4}$ &
			$2^{3}{:}\PSL_{2}(7)$ &
			&
			\cite{a:ABD-Un-CP}
			\\

			$3$ &
			$28$ &
			$36$ &
			$9$ &
			$7$ &
			$2$ &
			$\PSL_{2}(8)$ &
			$\D_{18}$ &
			$\D_{14}$ &
			$\PSL_{2}(8){:}3$ &
			&
			\cite{b:Handbook,a:Zhou-PSL2-non-sym}
			\\
			$4$ &
			$10$ &
			$15$ &
			$9$ &
			$6$ &
			$5$ &
			$\PSL_{2}(9)$&
			$3^2{:}4$&
			$\S_{4}$ &
			$\S_{6}$ &
			&
			\cite{a:Zhou-PSL2-non-sym,a:Zhou-CP-nonsym-alt}
			\\
			$5$ &
			$11$ &
			$11$ &
			$5$ &
			$5$ &
			$2$ &
			$\PSL_{2}(11)$ &
			$\PSL_{2}(5)$ &
			$\PSL_{2}(5)$ &
			$\PSL_{2}(11)$ &
			Hadamard &
			\cite{a:ABD-PSL2,a:Kantor-85-2-trans,a:Regueiro-reduction}
			\\
			$6$ &
			$12$ &
			$22$ &
			$11$ &
			$6$ &
			$5$ &
			$\M_{11}$ &
			$\PSL_{2}(11)$ &
			$\A_{6}$ &
			$\M_{11}$ &
			&
			\cite{a:Zhan-CP-nonsym-sprodic}
			\\
			$7$ &
			$22$ &
			$77$ &
			$21$ &
			$6$ &
			$5$ &
			$\M_{22}$ &
			$\PSL_{3}(4)$ &
			$2^{4}{:}\A_{6}$ &
			$\M_{22}$ &
			&
			\cite{a:Zhan-CP-nonsym-sprodic}
			\\
			$8$ &
			$22$ &
			$77$ &
			$21$ &
			$6$ &
			$5$ &
			$\M_{22}{:}2$ &
			$\PSL_{3}(4){:}2$ &
			$2^{4}{:}\S_{6}$ &
			$\M_{22}{:}2$ &
			&
			\cite{a:Zhan-CP-nonsym-sprodic}
			\\
			$9$ &
			$10$ &
			$15$ &
			$9$ &
			$6$ &
			$5$ &
			$\S_{6}$&
			$\S_{3}^{2}{:}2$&
			$\S_{4}{\times}2$ &
			$\S_6$ &
			&
			\cite{a:Zhou-PSL2-non-sym,a:Zhou-CP-nonsym-alt}
			\\
			$10$ &
			$15$ &
			$35$ &
			$7$ &
			$3$ &
			$1$ &
			$\A_7$  &
			$\PSL_2(7)$ &
			$(3\times \A_{4}){:}2$ &
			$\A_7$ &
			$\PG(3,2)$ &
			\cite{b:Handbook,a:Delandsheer-Lin-An-2001,a:Zhou-CP-nonsym-alt}
			\\
			$11$ &
			$15$ &
			$35$ &
			$7$ &
			$3$ &
			$1$ &
			$\A_{8}$ &
			$2^3{:}\PSL_3(2)$&
			$\A_{4}^{2}{:}2{:}2$ &
			$\A_{8}$ &
			$\PG(3,2)$&
			\cite{b:Handbook,a:Delandsheer-Lin-An-2001,a:Zhou-CP-nonsym-alt}
			\\
			\noalign{\smallskip}\hline\noalign{\smallskip}
			Note: & \multicolumn{11}{l}{The subgroups $G_{\alpha}$ and $G_{B}$ are point-stabilizer and block-stabilizer, respectively.  }\\
			
		\end{tabular}
	}
\end{table}

\section{Affine type}\label{sec:HA}

In this section, we provide some examples of $2$-designs with $\gcd(r,\lambda)=1$ admitting a flag-transitive automorphism group $G$ whose socle $T$ is an elementary abelian $p$-group of order $p^d$, $d\geq1$. All examples presented here can be read off from \cite{a:BDDKLS90} and therein references for $\lambda=1$, and in \cite{a:Biliotti-CP-nonsol-HA,a:Biliotti-CP-sym-affine,a:Biliotti-CP-sol-HA} for $\lambda>1$.\\
Since $G$ acts point-primitively on $\mathcal{D}$, the point set of $\mathcal{D}$ can be identified with $V_d(p)$, the $d$-dimensional vector space over the prime field $\GF(p)$ in a way that $T$ is the translation group of $V_d(p)$ and $G=TG_0$ is a subgroup of $\A\GL_d(p)$ with $G_0$  acting irreducibly on $V_d(p)$.
For each divisor $n$ of $d$ the group $\GaL_n(p^{d/n})$ has a natural irreducible action on $V_n(p^{d/n})$, thus we may choose $n$ to be minimal such that $G_{0}\leq \GaL_n(p^{d/n})$ in this action and write $q=p^{d/n}$.

\bigskip

To better understand some of the examples described below, some basics on $t$-spreads of vector spaces are provided. 
A (vectorial) $t\,$\emph{-spread} $\mathcal{S}$ of a $h$-dimensional vector space $V$ over the field $\GF(s)$, $s$ power of a prime, is a set of $t$-dimensional subspaces of $V$ partitioning $V-\left\{ 0\right\} $. Clearly, a $t\,$-spread of $V$ exists only if $t\mid h$. The incidence structure $\mathcal{A}(V,\mathcal{S})$, where the points are the vectors of $V$, the lines are the additive cosets of the elements of $\mathcal{S}$, and the incidence is the set-theoretic inclusion, is an \emph{Andr\'{e} translation structure} with associated translation group 
$T=\left\{ x\rightarrow x+w:w\in V\right\} $ and lines of size $s^{t}$. In particular $\mathcal{A}(V,\mathcal{S})$ is a $2$-$(s^{h},s^{t},1)$ design. Moreover, if $h$ is even and $t=h/2$, then $\mathcal{A}(V,\mathcal{S})$ is a \emph{translation plane} of order $s^{h/2}$.
If $\mathcal{S}$ and $\mathcal{S}^{\prime }$ are two $t$-spreads of $V$ such that $\psi $ is an
isomorphism from $\mathcal{A}(V,\mathcal{S})$ onto $\mathcal{A}(V,\mathcal{S}^{\prime })$ fixing the zero vector, then $\psi \in \GaL(V)$ and $\mathcal{S}^{\psi }=\mathcal{S}^{\prime }$. The converse is also true. Let $H_{0}\leq \GaL(V)$ preserving $\mathcal{S}$, then $TH_{0}$ is the automorphism group of $\mathcal{A}(V,\mathcal{S})$. The subset of $End(V,+)$ preserving each component of $\mathcal{S}$ is a field called the \emph{kernel} of $\mathcal{S}$ and is denoted by 
$\mathcal{K}(V,\mathcal{S})$, or simply by $\mathcal{K}$. Clearly, $\mathcal{K}-\left\{ 0\right\} \leq H_{0}$. Each component of $\mathcal{S}$ is a vector subspace of $V$ over $\mathcal{K}$ and $H_{0}\leq \GaL_{\mathcal{K}}(V)$. In particular, $\mathcal{A}(V,\mathcal{S})$ is a desarguesian affine space if, and only if, $\dim _{\mathcal{K}}Y=1$ for each $Y\in \mathcal{S}$. When this occurs the spread $\mathcal{S}$ is called \emph{regular}. More information on $t$-spreads and Andr\'{e} translation structures can be found in 
\cite{a:Andre-61,b:Dembowski,b:Johnson,b:Lueneburg}.\medskip

For $\lambda=1$, the examples of the designs and their automorphism groups are given in  Examples~\ref{ex:Desar-affine}-\ref{ex:Hering-affine} below:

\subsection{Desarguesian affine linear spaces}\label{ex:Desar-affine}
This design is the affine linear space $\AG_n(q)$ with $n\geq2$ and $q^n=p^d$, and $G_0\leq \GaL_n(q)$. If $Z$ is the centre of $\GL_n(q)$, then by the flag-transitivity of $G$, the group $ZG_0$ is transitive on the non-zero vectors of $\Dmc$, and so is one of the transitive linear
groups determined by Hering (see \cite[Appendix 1]{a:Liebeck-HA-rank3} for a list of these). From this, it follows quite easily that if $n>1$, then one of the following holds:
\begin{enumerate}[\rm (a)]
	\item $G$ is $2$-transitive on $V$ (hence given in \cite[Appendix 1]{a:Liebeck-HA-rank3});
	\item $n=2$, $q$ is $11$ or $23$, and $G$ is one of three soluble flag-transitive groups given in \cite[Table II]{a:Foulser-64};
	\item $n=2$, $q=9,11,19,29$ or $59$, $G_0^{(\infty)}=2.\A_5$ (where $G_0^{(\infty)}$ is the last term in the derived series of $G_0$), and $G$ is given in \cite[Table~II]{a:Foulser-64};
	\item $n=4$, $q=3$ and $G_0^{(\infty)}=2.\A_5$.
\end{enumerate}

\subsection{Non-Desarguesian translation affine planes}\label{ex:NonDesar-affine}
The examples for these designs with $\lambda=1$ are:
\begin{enumerate}[\rm (a)]
	\item The L\"{u}neburg planes constructed in \cite[Section 23]{b:Lueneburg}. These are affine planes of order $q^2$, where $q=2^{2a+1}\geq8$, and $^{2}\!\B_{2}(q)\vartriangleleft G_0\leq \Aut(^{2}\!\B_{2}(q))$;
	\item The Hering plane of order $27$ constructed in \cite{a:Hering-69}. Here $G_0=\SL_2(13)$ and $G$ is $2$-transitive on the points of $\Dmc$;
	\item The nearfield plane of order $9$. Here there are seven possibilities for $G$, given in \cite[Section 5]{a:Foulser-64}.
\end{enumerate}

\subsection{Hering spaces}\label{ex:Hering-affine}
These are two flag-transitive linear spaces with parameters $(3^6,3^2,1)$, constructed in \cite{a:Hering}. In both cases $G_0=\SL_2(13)$ and $G$ is $2$-transitive on the points. \medskip

For $\lambda>1$, the examples of the designs and their automorphism groups are given in  Examples~\ref{ex:affsub}-\ref{ex:tensprod} below:

\subsection{Designs where the blocks are subspaces of $\AG_{d}(p)$}\label{ex:affsub}
The incidence structure $\mathcal{D}=(V_{d}(p),\ell ^{G})$, where $\ell$ is a $u$-dimensional subspace of $V_{d}(p)$, is a $2$-design with parameters $\left( p^{d},p^{u},\frac{p^{u}-1}{p^{\gcd (u,d/n)}-1}\right)$  provided that ${\gcd (u,n,d/n)=1}$ and one of the following holds:

\begin{enumerate}[\rm (a)]
	
	\item $\gcd (u,d)<u<d/n$ and $\ell$ is contained in a 1-dimensional subspace of $V_{n}(q)$;
	
	\item $d-d/n\leq u<d$, $\gcd (u,d)<u$ and $\ell$ contains a hyperplane of $V_{n}(q)$;
	
	\item $(q,n,u)=(3,6,3)$ and $G_{0}\cong \SL_{2}(13)$.
	
\end{enumerate}

In (b), $\ell$ is the image under a polarity $\wp$ of $\PG_{d-1}(p)$ of either a block as in (a) or of a component of a $t$-spread $\mathcal{S}$ of $V_{d}(p)$ (the set of $1$-dimensional subspaces of $V_{n}(q)$ is a regular $d/n$-spread of $V_{d}(p)$). In (a) and (b), the possibilities for $G_{0}$ are as follows:

\begin{enumerate}[\quad \rm (i)]
	\item $\SL_{n}(q)\unlhd G_{0} \leq \GL_n(q) \cdot \langle \sigma_{0}\rangle$, where $\langle \sigma_{0}\rangle$ is the stabilizer in $\langle \sigma\rangle$ of $\ell$ and $\langle \sigma\rangle$ is induced by the field automorphism group of $\GF(q)$;
	
	\item $n$ is even and $Sp_{n}(q)\unlhd G_{0}\leq \GSp_n(q) \cdot \langle \sigma_{0}\rangle$;
	
	\item $q=2^{d/6}$, $n=6$ and $\G_{2}(q)^{\prime }\unlhd G_{0}\leq \G_{2}(q)^{\prime}\cdot \langle \sigma_{0}\rangle$;
	
	\item $(q,n,u)=(2,4,3)$ and $G_{0}\cong \A_{c}$, where $c=6,7$, or $G_{0}\cong \S_{6}$;
	
	\item $(q,n,u)=(3,4,3)$ and either $(\D_{8}\circ \Q_{8})\cdot 5\leq G_{0}\leq (\D_{8}\circ \Q_{8})\cdot \S_{5}$ or \\ $2\cdot \S_{5}^{-}\unlhd G_{0}\leq (2\cdot \S_{5}^{-}):2$, or $8\circ \SL_{2}(5)\unlhd G_{0}\leq (8\circ \SL_{2}(5)):2$;
	
	\item $(q,n,u)=(3,6,4)$, where $\ell$ is the image under $\wp $ of a component of any of the two $2$-spreads of $V_{6}(3)$ defining the two Hering spaces in \eqref{ex:Hering-affine}, and $G_{0}\cong \SL_{2}(13)$;
	
	\item $(q,n,u)=(3,6,5)$ and $G_{0}\cong \SL_{2}(13)$.
\end{enumerate}

\subsection{Designs where the blocks are union of distinct  subspaces of $\AG_{d}(p)$}\label{ex:unionaffsub}
These are four families of 2-$\left( q^n, p^u\omega, p^u\omega-1\right)$ designs $\mathcal{D}=(V_{d}(p),\ell ^{G})$, where $\ell$ is the union of $\omega$ cosets of a $u$-dimensional subspace $W$ of $V_{d}(p)$:

\begin{enumerate}[\rm (a)]
	\item 
	$\omega$ is a divisor of $p^{gcd(u,d/n)}-1$ such that $\gcd\left(\frac{p^{d}-1}{\omega},p^u\omega-1\right)=1$ and one of the following holds:
	
	\begin{enumerate}[(i)]
		\item
		$0\leq u<d/n$, $\A\GL_{1}(q) \unlhd G_X^X \leq \AGaL_{1}(q)$, where $X$ is any $1$-dimensional subspace of $V_{n}(q)$, and $\ell$ is any regular orbit of a Frobenius subgroup of $\A\GL_{1}(q)$ of order $p^{u}\omega$ contained in $X$; 
		\item
		$d-d/n\leq u<d$, there is a hyperplane $Y$ of $V_{n}(q)$ contained in $W$ and $G_Y$ acts on $V_n(q)/Y$ inducing a subgroup of $\AGaL_{1}(q)$ containing $\A\GL_{1}(q)$. The block $\ell$ is the union of $u-d+d/n$ cosets of $Y$ permuted regularly by a Frobenius subgroup of $\A\GL_{1}(q)$ of order $p^{u-d+d/n}\omega$. 
	\end{enumerate}
	
	In case (i), the translation complement $G_{0}$ is as in (i), (ii) and (iii) of \eqref{ex:affsub}, whereas in case (ii), $G_{0}$ is as in (i), (ii) and (vii) of \eqref{ex:affsub}.
	
	\item
	$p=2$, $n=6$, $u=d/3$, $\omega$ is any divisor of $q-1$ with $\gcd(q^{6}-1,q^{2}\omega -1)=1$. Here, $G_{2}(q)^{\prime }\unlhd G_{0}\leq G_{2}(q)^{\prime }\cdot \langle \sigma_{0}\rangle$, where $\langle \sigma_{0}\rangle$ is the stabilizer in $\langle \sigma\rangle$ of $\ell$ and $\langle \sigma\rangle$ is induced by the field automorphism group of $GF(q)$, and $W$ is any totally isotropic subplane of $V_6(q)$. The group $G_W \cap G_{2}(q)^{\prime } = [q^{5}]:\GL_2(q)$ preserves a non-isotropic $3$-dimensional subspace $Z$ containing $W$ inducing $\GL_{1}(q)$ on $Z/W$. The block $\ell$ is the union of $\omega$ cosets of $W$ in $Z$ permuted regularly by the subgroup of $\GL_{1}(q)$ of order $\omega$;
	
	\item 
	$\ell =V_{d}(p)-W$, where $W$ is a component of a $u$-spread $\mathcal{S}$ of $V_{d}(p)$, $\omega=p^{d-u}-1$ and one of the following occurs:
	
	\begin{enumerate}[(i)]
		\item $\mathcal{S}$ is a regular $d/2$-spread of $V_{d}(p)$, $\mathcal{D}$ is the complement of $\AG_{2}(q)$ and the following hold:
		\begin{enumerate}[(1)]
			\item $q=3,5,7,11,23$ and $(q-1)\cdot \A_{4}\unlhd G_{0}\leq {(q-1)}\cdot \S_{4}$;
			
			\item $q=9$ and either $2\cdot \S_{5}^{-}\unlhd G_{0}\leq (2\cdot \S_{5}^{-}):Z_{2}$ or $8\circ \SL_{2}(5)\unlhd G_{0}\leq (8\circ \SL_{2}(5)):2$;
			
			\item $q=11,19,29,59$ and $\SL_{2}(5)\unlhd G_{0}<{(q-1)}\circ \SL_{2}(5)$.
		\end{enumerate}
		\item $\mathcal{S}$ is a Hall $2$-spread of $V_{4}(3)$, $\mathcal{D}$ is the complement of the Hall plane of order $9$ and $G_{0}\cong 2\cdot \S_{5}^{-}$;
		\item $\mathcal{S}$ is the Hering $3$-spread of $V_{6}(3)$ and $\mathcal{D}$ is the $2$-$(3^{6},702,701)$ complement design of the Hering translation plane of order $3^{3}$ described in (b) of \eqref{ex:NonDesar-affine}, and $G_{0}\cong \SL_{2}(13)$.
	\end{enumerate}
	
	\item[(d)] The remaining cases:
	
	\begin{enumerate}[(i)]
		
		\item $\mathcal{D}=(V_{2}(7),\ell ^{G})$, where $\ell =W\cup (W+x)$ and $\dim W=1$, is a $2$-$(7^{2},14,13)$ design and $2\cdot \S_{4}^{-}\unlhd G_{0}\leq 6\cdot \S_{4}^{-}$. Here, $\ell ^{G}$ is the set of all pairs of parallel lines of $\AG_{2}(7)$;
		
		\item $\mathcal{D}=(V_{4}(3),\ell ^{G})$, where $\ell =W\cup (W+x)$ and $\dim W=2$, is a $2$-$(3^{4},18,17)$ design, and the following hold:
		
		\begin{enumerate}[(1)]
			\item $G_{0}\cong 2\cdot \S_{5}^{-}$ and $\ell ^{G}$ is the set of all pairs of parallel lines of the Hall plane of order $9$;
			
			\item $2\cdot \S_{5}^{-}\unlhd G_{0}\leq (2\cdot \S_{5}^{-}):2$ or $8\circ \SL_{2}(5)\unlhd G_{0}\leq (8\circ \SL_{2}(5)):2$ and $\ell ^{G}$ is either the set of all pairs of parallel lines of $\AG_{2}(9)$, or $\left\vert W^{G_{0}}\right\vert =40$ and both $W$ and $W+x$ are two subplanes of $\AG_{4}(3)$ preserved by a Sylow $3$-subgroup of $G_{0}$.
		\end{enumerate}
		
		\item $\mathcal{D}=(V_{6}(3),\ell ^{G})$ is a $2$-$(3^{6},2\cdot 3^{u},2\cdot 3^{u}-1)$ design, where $u=1,2,3$, and $G_{0}\cong \SL_{2}(13)$. Here, $\ell =W\cup (W+x)$ and $W$ is any $u$-dimensional subspace of $V_{6}(3)$ such that either $\left\vert W^{G_{0}}\right\vert =364$, or $u=2$ and $\left\vert W^{G_{0}}\right\vert =91$, or $u=3$ and $\left\vert W^{G_{0}}\right\vert =28$. The number of isomorphism classes of such designs is $4u$ for $u=1,2,3$ and $\left\vert W^{G_{0}}\right\vert =364$, five for $u=2$ and $\left\vert W^{G_{0}}\right\vert =91$, and one for $u=3$ and $\left\vert W^{G_{0}}\right\vert =28$;
		
		\item $\mathcal{D}=(V_{6}(3),\ell ^{G})$ is a $2$-$(3^{6},324,323)$ design and $G_{0}\cong \SL_{2}(13)$. Let $W$ be any $4$-dimensional subspace of $V_{6}(3)$ such that $G_{0,W}\cong \Q_{8}:3$. Then the block $\ell =\cup _{i=1}^{4}(W+x_{i})$ where the $(W+x_{i})$'s are four $4$-dimensional subspaces of $\AG_{6}(3)$ preserved by a Sylow $3$-subgroup of $G_{0}$ and permuted transitively by a Klein subgroup of $G$. The number of isomorphism classes of such designs is three.
	\end{enumerate}
\end{enumerate}

\subsection{Designs where a block base is properly contained in a component of a transitive $t$-spread of $V_{d}(p)$} \label{ex:1comspread}

The group $G_{0}$ acts transitively on a $t$-spread $\mathcal{S}$ of $V_{d}(p)$ and $\mathcal{D}=(V_{d}(p),\ell ^{G})$, where $\ell \subseteq W-\{0\}$, $W\in \mathcal{S}$. More precisely, one of the following holds:

\begin{enumerate}[\rm (a)]
	
	\item $\ell $ is a quadrangle contained in the union of two $1$-dimensional subspaces of $W$. Hence, $\mathcal{D}$ is a $2$-$(q^{n},4,3)$ design, provided that one of the following holds:
	
	\begin{enumerate}[(i)]
		
		\item $(q,n)=(3,4)$, $\mathcal{S}$ is a Hall $2$-spread and either 
		$(\D_{8}\circ \Q_{8})\cdot 5\leq G_{0}\leq (\D_{8}\circ \Q_{8})\cdot \S_{5}$ 
		or $G_{0}\cong 2\cdot \S_{5}^{-}$;
		
		\item $(q,n)=(9,2)$, $\mathcal{S}$ is a regular $2$-spread of $V_{4}(3)$ and either $2\cdot \S_{5}^{-}\unlhd G_{0}\leq (2\cdot \S_{5}^{-}):2$ or $8\circ \SL_{2}(5)\unlhd G_{0}\leq (8\circ \SL_{2}(5)):2$;
		
		\item $(q,n)=(3,6)$, $\mathcal{S}$ is one of the two $2$-spreads defining the two Hering spaces in \eqref{ex:Hering-affine}, and $G_{0}\cong \SL_{2}(13)$.
	\end{enumerate}
	
	In cases (i),(ii) and (iii) the quadrangle $\ell$ is an orbit under a cyclic subgroup of order $4$ of $G_{W}\cong \SL_{2}(3)$. There are six of such quadrangles which are contained in $W-\{0\}$. Moreover, they are permuted transitively by $G_{W}$. Therefore, the $2$-designs in (i), (ii) and (iii) are unique up to isomorphism, since $G_{0}$ acts transitively on the $2$-spread $\mathcal{S}$.
	
	\item $\ell =W-\{0\}$, $\mathcal{D}$ is a $2$-$(p^{d},p^{t}-1,p^{t}-2)$ design and the following hold:
	
	\begin{enumerate}[(i)]
		
		\item $(q,n,t)=(3,4,2)$, $\mathcal{S}$ is a Hall $2$-spread and either $(\D_{8}\circ \Q_{8})\cdot 5\unlhd G_{0}\leq (\D_{8}\circ \Q_{8})\cdot \S_5$ or $G_{0}\cong 2\cdot \S_{5}^{-}$;
		
		\item $(q,n,t)=(3,4,2)$, $\mathcal{S}$ is a regular $2$-spread and either $2\cdot \S_{5}^{-}\unlhd G_{0}\leq (2\cdot \S_{5}^{-}):2$ or $8\circ \SL_{2}(5)\unlhd G_{0}\leq (8\circ \SL_{2}(5)):2$;
		
		\item $(q,n,t)=(3,6,3)$, $\mathcal{S}$ is a one of the two $2$-spreads defining the Hering spaces in \eqref{ex:Hering-affine}, and $G_{0}\cong \SL_{2}(13)$.
	\end{enumerate}
\end{enumerate}

\subsection{Designs where a block base is contained in at least two components of a transitive $t$-spread of $V_{d}(p)$} \label{ex:2comspread}
There are two examples in this case:
\begin{enumerate}[\rm (a)]
	\item $\mathcal{D}=(V_{4}(3),\ell ^{G})$, where $\ell $ is any of the two $\SL_{2}(5)$-orbits on the set of non-zero vectors, is a $2$-$(3^{4},40,39)$ design and $2\cdot \S_{5}^{-}\unlhd G_{0}\leq (2\cdot \S_{5}^{-}):2$ or $8\circ \SL_{2}(5)\unlhd G_{0}\leq (8\circ \SL_{2}(5)):2$. Here 
	\[
	\ell=\bigcup\limits_{i=1}^{5}(W_{i}-\{0\})=\bigcup\limits_{j=1}^{5}(W_{j}^{\prime }-\{0\}),
	\]
	where the $W_{i}$'s and the $W_{j}$'s are components of two distinct Hall 2-spreads $\mathcal{S}$ and $\mathcal{S}^{\prime }$ of $V_{4}(3)$. Moreover, $\mathcal{D}$ is unique up to isomorphism;
	
	\item $\mathcal{D}=(V_{6}(3),\ell ^{G})$, where $G_{0}\cong \SL_{2}(13)$ and $\ell$ is any orbit under a cyclic subgroup of order $26$, is a $2$-$(3^{6},26,25)$ design. There are eight pairwise non-isomorphic types of such $2$-designs;
	
	\item $\mathcal{D}=(V_{6}(3),\ell ^{G})$, where $G_{0}\cong \SL_{2}(13)$ and $\ell$ is any of the two orbits of length $52$ under a subgroup of $G_{0}$ isomorphic to $2\cdot(13:6)$, is a $2$-$(3^{6},52,51)$ design. There are two non-isomorphic $2$-designs. 
	Indeed, if $\mathcal{S}$ is the $3$-spread  of $V_{6}(3)$   yielding the Hering translation plane of order $3^{3}$, in one design $\ell=(W_{1} \cup W_{2})-\{0\}$, where $W_{1},W_{2} \in \mathcal{S}$, whereas in the other $\ell\cap(W-\{0\})=\{\pm x_{W}\}$ for each $W\in \mathcal{S}-\{W_{1},W_{2}\}$.
\end{enumerate}

\subsection{Designs arising from suitable subsets of $\AG_{2}(p)$ of cardinality prime to $p$}\label{ex:subbidim}

These are $2$-$(p^{2},\lambda +1,\lambda )$ designs $\mathcal{D}=(V_{2}(p),\ell ^{G})$, where $\gcd(p,\lambda +1)$=1, and the possibilities are as follows:

\begin{enumerate}[(a)]
	
	\item $p=5$, $\lambda =5,7$ and $\SU_{2}(3)\unlhd G_{0}\leq \GU_{2}(3)$;
	
	\item $p=7$, $\lambda =5,7,11,23$ and $2\cdot \S_{4}^{-}\unlhd G_{0}\leq 6\cdot \S_{4}^{-}$;
	
	\item $p=11$, and one of the following holds:
	\begin{enumerate}[(i)]
		\item $\lambda =7$ and either $5\times \SL_{2}(3)\unlhd G_{0}\leq 5\times \GL_{2}(3)$ or $G_{0}\cong \SL_{2}(5)$;
		\item $\lambda =11$ and $G_{0}\cong \SL_{2}(5)$;
		\item $\lambda =19$ and either $5\times \SL_{2}(3)\unlhd G_{0}\leq 5\times \GL_{2}(3)$ or $\SL_{2}(5)\unlhd G_{0}\leq 5\times \SL_{2}(5)$;
		\item $\lambda =23$ and $G_{0}\cong \SL_{2}(5)$;
		\item $\lambda =29$ and $5\times \SL_{2}(3)\unlhd G_{0}\leq 5\times \GL_{2}(3)$. 
	\end{enumerate}
	\item $p=19$, $\lambda =7,11,23,71,119$ and $G_{0}\cong 9\times \SL_{2}(5)$;
	
	\item $p=23$, $\lambda =5,7,23,43,47,65,175,263$ and $G_{0}\cong {11}\times 2\cdot \S_{4}^{-}$;
	
	\item $p=29$, $\lambda =11,13,23,41,83,139,167$ and $7\times \SL_{2}(5)\unlhd G_{0}\leq 28\circ \SL_{2}(5)$;
	
	\item $p=59$, $\lambda =7,11,19,23,119,173,289,347$ and $G_{0}\cong 29\times \SL_{2}(5)$.
\end{enumerate}

A detailed description of these cases is contained in \cite{a:Biliotti-CP-nonsol-HA,a:Biliotti-CP-sol-HA}.

\subsection{Example arising from tensor product decomposition of the vector space}\label{ex:tensprod}
In this case $G_{0}\cong \S_{3}\times \PSL_{3}(2)$ preserves the tensor decomposition $V_{6}(2)=V_{2}(2)\otimes V_{3}(2)$. Hence, $\mathcal{D}=(V_{6}(2),\ell ^{G})$, where $\ell =u\otimes \PG_{2}(2)$, $u\neq 0$ is a copy of a Fano plane, is a $2$-$(2^{6},7,2)$ design admitting $G=TG_{0}$ as a flag-transitive, point-primitive automorphism group. 

In each case, except for \eqref{ex:tensprod}, $G$ acts point-$2$-transitively on $\mathcal{D}$. All examples given here do occur and the construction details for each $2$-design can be read off from \cite{a:Biliotti-CP-nonsol-HA,a:Biliotti-CP-sol-HA}.

The previous paragraphs contain an exhaustive list of $2$-designs admitting a flag-transitive, point-primitive automorphism group except for the semilinear $1$-dimensional case. Some of the known examples of $2$-designs arising from the latter case are described in the following final paragraph. The list of examples presented below is far from being exhaustive.

\section{The semilinear $1$-dimensional case}\label{sec:onesemdim}

As noted in the Introduction, we have so far been unable to handle completely the case where $G\leq \AGaL_{1}(p^{d})$. This brief section is devoted to illustrate some examples, where $\mathcal{D}$ is a $2$-design with $p^{d}$ points and $\gcd(r,\lambda)=1$ admitting a subgroup $G$ of the semilinear $1$-dimensional group $\AGaL_{1}(p^{d})$ as a flag-transitive and point-primitive automorphism group. For $\lambda=1$, remarkable examples are due to Kantor and to Pauley and Bamberg. The reader is referred to \cite{a:BDDKLS90} and to \cite{a:Bamberg-08}, respectively, for more details on these constructions.

Suppose that $\lambda>1$. Case (a) in \eqref{ex:affsub} or in \eqref{ex:unionaffsub} 
occurs also for $n=1$ and $G\cong \AGL_{1}(p^{d}): \langle \sigma_{0}\rangle$, where $\langle \sigma_{0}\rangle$ is the stabilizer in $\langle \sigma\rangle$ of $\ell$, 
and $\langle \sigma\rangle$ is induced by the field automorphism group of $\GF(p^{d})$. 
Moreover, $\mathcal{D}$ as in \eqref{ex:tensprod} admits $G=TG_{0}$ and $21\trianglelefteq G_{0}\leq 21:3$ inside $\GaL_{1}(2^{6})$ as a flag-transitive point primitive automorphism group. Further remarkable examples of $2$-$(2^{6},7,2)$ designs involving the tensor decomposition of $V_{6}(2)=V_{2}(2)\otimes V_{3}(2)$ and admitting $G=TG_{0}$ and $21\trianglelefteq G_{0}\leq 21:3$ as a flag-transitive, point-primitive automorphism group are the following: $\mathcal{D}=(V_{6}(2),\ell_{h} ^{G})$, with 
$\ell_{h} =\{u_{1}\otimes w^{\gamma^{i}}+ u_{2}\otimes w^{\gamma^{i+h}}\mid i=1,\ldots,7\}$, 
where $\langle \gamma\rangle$ is cyclic of order $7$, $u_{1}\neq u_{2}$ and $w \neq 0$, for each $h=1,...,6$. There are not corresponding examples of such $2$-designs with $\gcd(r,\lambda)=1$ admitting a non-solvable automorphism group. More details on these $2$-designs are provided in \cite{a:Biliotti-CP-sym-affine}.

If $\mathcal{D}$ is symmetric and admits a flag-transitive and point-primitive automorphism group $G$ of affine type, then $G\leq \AGaL_{1}(p^{d})$ by \cite{a:Biliotti-CP-sym-affine}. Moreover, if $O(G)$ denotes the maximal subgroup of $G$ of odd order, then $O(T)=T:\left\langle \bar{\omega}^{i},\sigma ^{y}\bar{\omega}^{j}\right\rangle $, where $y\mid d$, $i\mid (p^{d}-1,j\frac{p^{d}-1}{%
	p^{y}-1})$ and $\bar{\omega}:x\rightarrow \omega x$, $\omega $ is a
generator of $\GF(p^{d})^{\ast }$, and $\sigma :x\rightarrow x^{p}$, acts flag-transitively, point-primitively on $\mathcal{D}$. Finally, $\mathcal{D}$ has parameters $\left( p^{d},\theta \frac{p^{d}-1}{i},\theta^{2}\frac{p^{d}-1-i/\theta}{i^{2}}\right)$, where $\theta$ is a divisor of $d/y$. Details on the structure of $\mathcal{D}$ and on the corresponding examples due to Paley, Chowla and Lehmer, are available in \cite{a:Biliotti-CP-sym-affine}.

\section*{Statements and Declarations}

The authors confirm that this manuscript has not been published elsewhere. It is not also under consideration by another journal. They also confirm that there are no known conflicts of interest associated with this publication.  The authors have no competing interests to declare that are relevant to the content of this article and they confirm that availability of data and material is not applicable. 
Shenglin Zhou
is supported by the NNSF of China (Grant No.11871224).



\end{document}